\documentclass[a4,psfig,english,10pt,epsf,portrait]{article}
\usepackage{amsmath,amsfonts,latexsym,amscd,amssymb,theorem}
\usepackage{epsfig}
\usepackage{amsmath}
\usepackage{psfrag}
\def\R{\hbox{\bf R}}
\def\Z{\hbox{\bf Z}}

\def\N{\hbox{\bf N}}

\def\e{\varepsilon}

\def\<{\langle}
\def\>{\rangle}
\newcommand{\ba}{\begin{eqnarray}}
\newcommand{\ea}{\end{eqnarray}}

\newtheorem{thm}{Theorem}[section]

\newtheorem{theorem}[thm]{Theorem}

\newtheorem{rem}[thm]{Remark}

\numberwithin{equation}{section}

\renewcommand{\R}{{\mathbb R}}
\renewcommand{\Z}{{\mathbb Z}}
\renewcommand{\N}{{\mathbb N}}

\begin{document}

\title{\bf Dislocation dynamics: from microscopic models to macroscopic crystal plasticity}
\author{
\normalsize\textsc{ A. El Hajj $^{*}$, H. Ibrahim $^{*,\ddagger,\dagger}$, R. Monneau
  \footnote{Universit\'{e} Paris-Est, CERMICS,
Ecole des Ponts, 6 et 8 avenue Blaise Pascal, Cit\'e
Descartes Champs-sur-Marne, 77455 Marne-la-Vall\'ee Cedex 2, France.
E-mails: ibrahim@cermics.enpc.fr, monneau@cermics.enpc.fr
\newline \indent $\,\,{}^\dagger$CEREMADE, Universit\'{e}
Paris-Dauphine, Place De Lattre de Tassigny, 75775 Paris Cedex 16,
France
\newline \indent $\,\,{}^\ddagger$LaMA-Liban,
Lebanese University, P.O. Box 826 Tripoli, Lebanon}}}
\vspace{20pt}

\maketitle

%%%%%%%%%%%%%%%%%%%%%%%%%%%%%%%%%%%%%%%%%%%%%%%%%%%%%%%%%%%%%%%%%%%%%%%%%%
%%%%%%%%%%%%%%%%%%%%%%%%%%%%%%%%%%%%%%%%%%%%%%%%%%%%%%%%%%%%%%%%%%%%%%%%%%

\centerline{\small{\bf{Abstract}}} \noindent{\small{In this paper we
    study the connection between four models describing dislocation dynamics: a generalized 2D
 Frenkel-Kontorova model at the atomic level, the Peierls-Nabarro model,
 the discrete dislocation dynamics and a macroscopic model with
 dislocation densities. We show how each model can be deduced from the previous one at a smaller scale.}}

\hfill\break
 \noindent{\small{\bf{AMS subject classifications:}}} {\small{35B27,
     35F20, 45K05, 47G20, 49L25}}\hfill\break 
  \noindent{\small{\bf{Key words:}}} {\small{discrete to continuum
      approach, Frenkel-Kontorova model, phase field model, particle
      systems, periodic homogenization, Hamilton-Jacobi equations, non local equations}}\hfill\break

\section{Introduction}

In this paper, we focus on the modelling of dislocation dynamics. 
We refer the reader to the book of Hirth and Lothe \cite{HL} for a detailed introduction to dislocations.
Our study ranges from atomic models to macroscopic crystal plasticity.
At each scale, dislocations can be described by a suitable model.
Our goal is to explain how we can deduce a model at a larger scale, 
from the model at a smaller scale.

Even if our derivation will be done on some simplified models (essentially 2D and 1D models),
we hope that our contribution will shed light, even on some well-known models.
More precisely, we will consider the following four models, from the smaller to the larger scale:\\
1. Generalized Frenkel-Kontorova  model (FK)\\
2. Peierls-Nabarro model (PN)\\
3. Dynamics of discrete dislocations (DDD)\\
4. Dislocation density model (DD)

Schematically, the four models are related as shown below (see also
Figure 6.1 for a more detailed diagram):

 \begin{equation}\label{tableau}{\fbox{(FK)$_{\e_{1}, \e_{2},
        \e_{3} >0}$ }  \overset{\e_1 \to 0}{   \longrightarrow  }
\fbox{(PN)$_{\e_{2},
        \e_{3} >0}$}    \overset{\e_2 \to 0}{ \longrightarrow }
\fbox{(DDD)$_{\e_{3} >0}$}  \overset{\e_3 \to 0}{ \longrightarrow }
\fbox{(DD)}}\end{equation}  

The rest of the paper is composed of four sections.
Each section presents one model, and explains how this model 
can be deduced from the previous model at a smaller scale.

%%%%%%%%%%%%%%%%%%%%%%%%%%%%%%%%%%%%%%%%%%%%%%%%%%%%%%%%%%%%%%%%%%%%%%%%%%
%%%%%%%%%%%%%%%%%%%%%%%%%%%%%%%%%%%%%%%%%%%%%%%%%%%%%%%%%%%%%%%%%%%%%%%%%%
\section{Generalised Frenkel-Kontorova model}\label{s1}

%%%%%%%%%%%%%%%%%%%%%%%%%%%%%%%%%%%%%%%%%%%%%%%%%%%%%%%%%%%%%%%%%%%%%%%%%%
\subsection{Geometrical description}

We call  $(e_1,e_2,e_3)$ a direct orthonormal basis of the threedimensional space.
We consider a perfect crystal $\Z^3$ where each position with integer coordinates is occupied by one atom.
We want to describe dislocations, which are certain ``line defects'' in the crystal. To simplify the presentation, 
we will assume that the material is invariant by integer translations in the direction $e_3$. 
Because of this assumption, we can simply consider the cross section of the crystal in the plane $(e_1,e_2)$
where each atom is now assumed to have a position $I\in\Z^2$ in the perfect crystal.
We also assume that each atom $I$ can have a displacement $U_I \in\R$ in the direction $e_1$, 
such that the effective position of the atom $I$ is $I+ U_I e_1$. 

%%

%  \begin{figure}[h]
%  \centering\epsfig{figure=F1.eps,width=60mm}
%  \caption{Perfect crystal}\label{F1}
%  \end{figure}

%  \begin{figure}[h]
%  \centering\epsfig{figure=F2.eps,width=60mm}
%  \caption{Schematic view of a edge dislocation in the crystal}\label{F2}
%  \end{figure}

\begin{figure}[!h]
    \begin{minipage}[b]{.46\linewidth}
     \centering\epsfig{figure=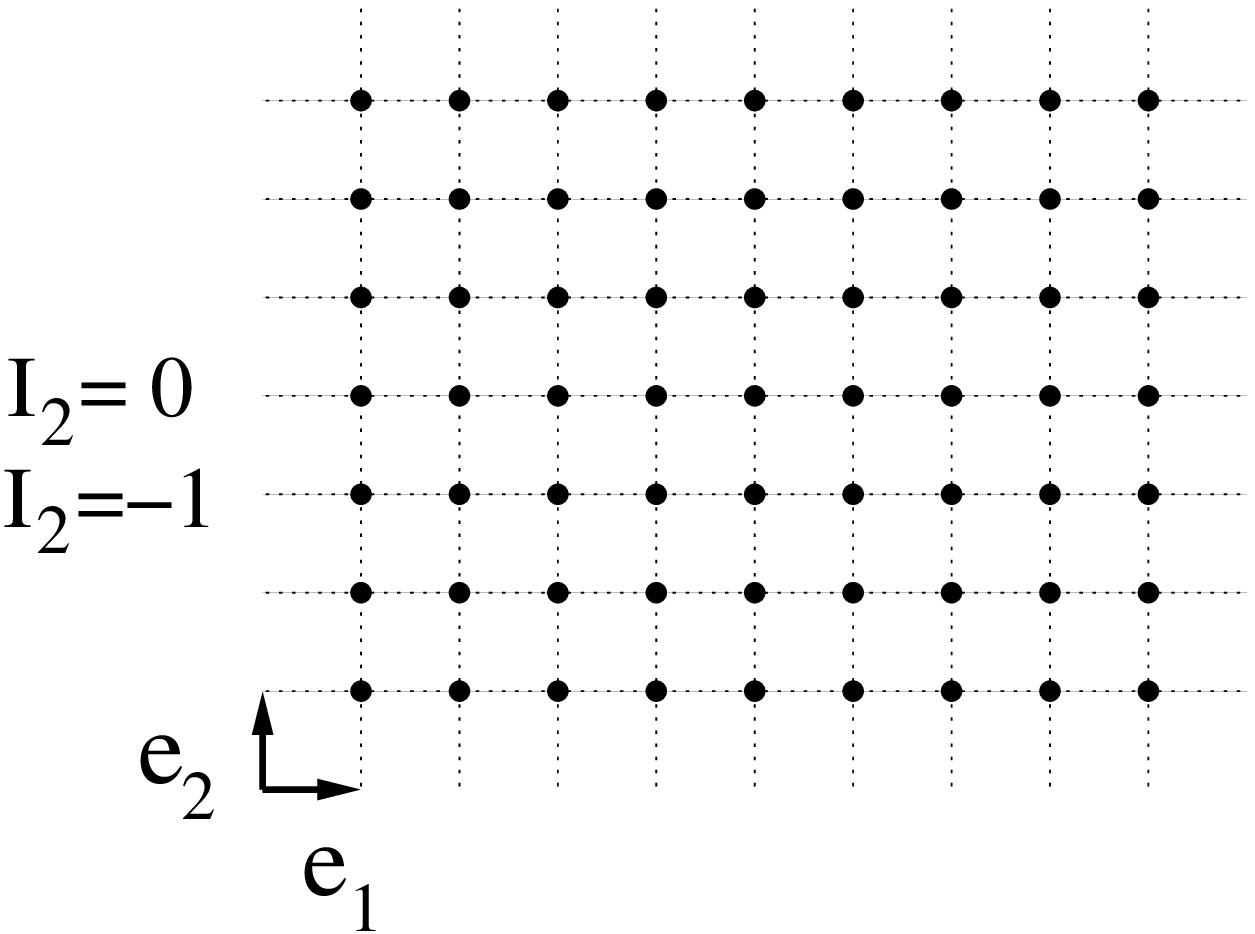,width=\linewidth}
     \caption{Perfect crystal 
      $\quad \quad \quad \quad \quad \quad \quad \quad \quad$ \label{F1}}
    \end{minipage} \hfill
    \begin{minipage}[b]{.46\linewidth}
     \centering\epsfig{figure=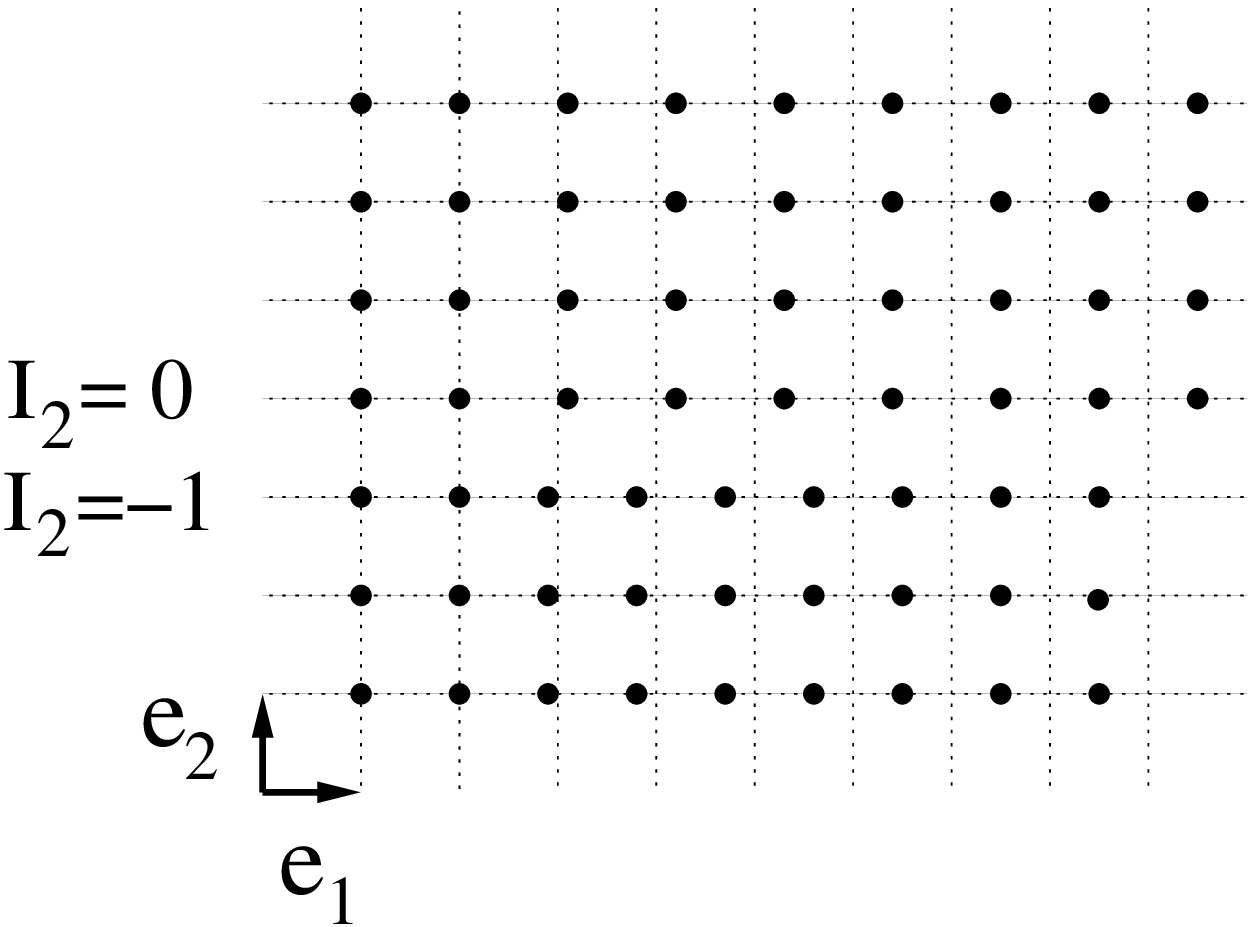 ,width=\linewidth}
     \caption{Schematic view
  of  a edge dislocation in the crystal\label{F2}}
    \end{minipage}
   \end{figure}

On Figure \ref{F1} below is represented a view of the perfect crystal. 
On Figure \ref{F2} we can see a schematic view of a edge dislocation in the crystal.
On this picture, the upper part $\left\{I_2\ge 0\right\}$ of the crystal has been expanded 
to the right of a vector  $\frac12 e_1 $, 
while the lower part $\left\{I_2\le  -1\right\}$ of the crystal has been contracted to the left
of a vector $-\frac12 e_1$. The net difference between these two vectors is $e_1$ and is called the Burgers vector of this dislocation.

In order to describe a edge dislocation in our formalism, let us make a few assumptions.
We will assume that the dislocation defects are essentially described by the mismatch between the two planes  
$I_2=0$ and $I_2=-1$, like on  Figure \ref{F2}. For this reason, and also in order to simplify the analysis, 
we assume that the displacement of the crystal satisfies the following {\it antisymmetry property}
\begin{equation}\label{eq::40}
U_{(I_1,-I_2)}=-U_{(I_1,I_2-1)} \quad \mbox{for all}\quad I=(I_1,I_2)\in\Z^2.
\end{equation}
Let us also define the discrete gradient
$$(\nabla^d U)_I=\left(\begin{array}{c}
U_{I+e_1}-U_I\\
U_{I+e_2}-U_I
\end{array}\right).$$
Remark that defects in the crystal can be seen as regions where the discrete gradient is not small.\\

\noindent {\bf Formalism for a edge dislocation with Burgers vector $e_1$}\\
In our formalism, a edge dislocation like the one of Figure \ref{F2}, can be represented
by a displacement $U_I$ satisfying
$$\left\{\begin{array}{l}
U_{(I_1,0)}=-U_{(I_1,-1)} \to 0 \quad \mbox{as}\quad I_1 \to
-\infty\\
\\
\displaystyle{U_{(I_1,0)}=-U_{(I_1,-1)} \to \frac12  \quad \mbox{as}\quad
I_1 \to +\infty} .
\end{array}\right.$$
Because we assume that the dislocation core lies in the two planes $I_2=0$ and $I_2=-1$, 
it is reasonable to assume that 
all the components of the discrete gradient are small, 
except components $U_{I+e_2}-U_I$ for $I=(I_1,I_2)$ with $I_2=-1$. More precisely, 
we assume that there exists a small $\delta >0$ such that 
\begin{equation}\label{eq::1}
\left\{\begin{array}{l}
|U_{I+e_1}-U_I|\le \delta \quad \mbox{for all}\quad  I=(I_1,I_2) \in\Z^2     \\
|U_{I+e_2}-U_I|\le \delta \quad \mbox{for all}\quad  I=(I_1,I_2) \in\Z^2\quad \mbox{with}\quad I_2\not= -1 .
\end{array}\right.
\end{equation}
Moreover, if there is no applied stress on the crystal, then it is
reasonable to assume that
$$\mbox{dist }\left((\nabla^d U)_I, \Z^2 \right) \to 0 \quad \mbox{as}\quad
|I|\to +\infty .$$

%%%%%%%%%%%%%%%%%%%%%%%%%%%%%%%%%%%%%%%%%%%%%%%%%%%%%%%%%%%%%%%%%%%%%%%%%%
\subsection{The energy and the dynamics}

We assume that the energy of a configuration $U=(U_I)_{I\in\Z}$ of the crystal can be formally written as
$$E(U)=\frac12 \sum_{I\not= J}\widetilde{W}(U_I-U_J)$$
where $\widetilde{W}:\R\to \R$ is a potential describing nearest neighbors interactions satisfying\\
\noindent {\bf Assumption ($\widetilde{A1}$)}
$$\left\{\begin{array}{lll}
\mbox{\bf i) (Regularity)}&\quad \widetilde{W}\in C^3(\R)\\
\mbox{\bf ii) (Periodicity)}&\quad \widetilde{W}(a+1)=\widetilde{W}(a) &\quad \mbox{for all}\quad a\in\R\\
\mbox{\bf iii) (Minimum on $\Z$)}&\quad \widetilde{W}(\Z)=0
<\widetilde{W}(a)  &\quad \mbox{for all}\quad
 a\in\R\backslash \Z\\
\mbox{\bf iv) (Local harmonicity of $\widetilde{W}$)}&\quad
\widetilde{W}(a)=\frac12 a^2 &\quad \mbox{for all}\quad
 |a|< \delta
\end{array}\right.$$
where $\delta >0$ is introduced in (\ref{eq::1}).
Remark that the periodicity of the potential $\widetilde{W}$ reflects the periodicity of the crystal, 
while the mimimum property of $\widetilde{W}$ is consistent with the
fact that the perfect crystal $\Z^2$ 
is assumed to minimize its energy. Assumption iv) will be used for later simplification.\\

Then we assume that we are in a regime where the crystal reaches very quickly the equilibrium 
in the regions where there is no defects,  
i.e. satisfies 
%(in particular with the escape at infintity of high speed elastic waves)
\begin{equation}\label{eq:1}
0=-\nabla_{U_I} E(U) \quad \mbox{for all}\quad I=(I_1,I_2)\in\Z^2\quad \mbox{with}\quad I_2\not= 0,-1
\end{equation}
while we have the following fully overdamped dynamics in the two planes where the dislocation lives
(describing the average friction of the lattice on the effective dissipative motion of the dislocations):
\begin{equation}\label{eq:2}
\frac{d}{dt} U_I =-\nabla_{U_I} E(U) \quad \mbox{for all}\quad I=(I_1,I_2)\in\Z^2\quad \mbox{with}\quad I_2 = 0,-1.
\end{equation}
Let us mention that we do not have a fundamental justification of this dynamics, 
but we think that one of the main justification of this model is that other known models at larger scales 
can be deduced from this particular model. For physical justifications
of the dissipative effects in the motion of dislocations, see
\cite{AI,HL}. See also \cite{KT1,KT,KT2} for a fundamental justification of the overdamped dynamics based on explicit computations in a 1D Hamiltonian model.

Taking into account the local harmonic assumption ($\widetilde{A1}$) iv), 
applied where the components of the discrete gradient are small (see (\ref{eq::1})), joint to the antisymmetry property
defined in (\ref{eq::40}), 
we can rewrite system (\ref{eq:1})-(\ref{eq:2}) as follows for all $t>0$:
\begin{equation}\label{eq:3}
\left\{\begin{array}{ll}
\displaystyle{0=\sum_{J\in\Z^2, \ |J-I|=1} (U_J-U_I)} & \quad \mbox{for all}\quad I=(I_1,I_2)\in\Z^2 \quad \mbox{with}\quad I_2\ge 1\\
\\
\displaystyle{\frac{d}{dt} U_I=-\widetilde{W}'(2U_I) + \sum_{J\in\Z^2,\ |J-I|=1,\  J_2\ge 0} (U_J-U_I)} &\quad \mbox{for all}\quad I=(I_1,I_2)\in\Z^2\quad \mbox{with}\quad I_2 = 0 .
\end{array}\right.
\end{equation}

We call this model a generalised Frenkel-Kontorova model. Even if this system of equations is not standard, it is nevertheless possible to define a unique solution under suitable assumptions, in the framework of viscosity solutions (see \cite{FinoIM}).
We refer the reader to the book of Braun, Kivshar, \cite{BK} for a detailed presentation of the classical FK model.
For homogenization results of FK models, we refer the reader to
\cite{FIM2}. For the
description of vacancy defects at equilibruim, see \cite{HMR2007}. See
also \cite{HMR2006}, where the authors study the problem involving a
dislocation inside the interphase between two identical lattices. Their
model corresponds to our model (\ref{eq:3}) at the equilibrium with the
potential $\widetilde{W}$ is a cosinus function. For other 2D FK models, see \cite{CB,CCHO}.

\begin{rem}
It is important to remark that we used condition (\ref{eq::1}) only to derive the model.
We do not know and we do no claim that there exists solutions of system (\ref{eq:3}) satisfying condition (\ref{eq::1}).
From now on, we only consider solutions of system (\ref{eq:3}) without requiring further assumptions on the solutions.   
\end{rem}

\begin{rem}
When we freeze the components $U_I=0$ for $I_2\ge 1$, and change the evolution equation forgetting the index $J$ with $J_2=1$, this leads to the following classical fully overdamped Frenkel-Kontorova model satisfied by $V_i:=U_{(i,0)}$
$$\frac{d}{dt} V_i=V_{i+1}+V_{i-1} -2V_i-\widetilde{W}'(2V_i) .$$   
\end{rem}

%%%%%%%%%%%%%%%%%%%%%%%%%%%%%%%%%%%%%%%%%%%%%%%%%%%%%%%%%%%%%%%%%%%%%%%%%%
\subsection{The asymptotic stress created by a single dislocation}

In this subsection, we will compute the asymptotic stress created by a single dislocation.
To this end, we first compute the effective Hook's law of the lattice.

\noindent {\bf Computation of the Hook's law}\\
Let us consider an affine displacement
$$U_I=a\cdot I + C \quad \mbox{with}\quad a=(a_1,a_2)\in\R^2$$
where $C\in\R$ is a constant.
Then the energy by unit cell is
$$\displaystyle{{\mathcal E}=\widetilde{W}(U_{I+e_1}-U_I))+\widetilde{W}(U_{I+e_2}-U_I))=\frac12 (a_1^2 + a_2^2)}$$
for $|a|< \delta$.
Reminding the fact that $U$ is the displacement in the $e_1$ direction, we get that the strain $e$ (i.e. the symmetric part of the gradient of the displacement) is given by
$$e=  \left(\begin{array}{ll}
e_{11} & e_{12}\\
e_{21} & e_{22}
\end{array}\right)= \frac12 \left( \nabla U \otimes e_1 + e_1 \otimes\nabla U \right) =\left(\begin{array}{ll}
a_1   &   a_2/2\\
a_2/2 &    0
\end{array}\right) .$$
Therefore
$$\displaystyle{{\mathcal E}(e)=\frac12 e_{11}^2 + 2 e_{12}^2} .$$
Recalling that the stress is given by $\displaystyle{\sigma^0 = 
\frac{\partial {\mathcal E}}{\partial e}}$, we get the Hook's law:
$$\sigma^0=\left(\begin{array}{ll}
e_{11} & 2e_{12}\\
2e_{21} & 0
\end{array}\right) .$$

\noindent {\bf Computation of the stress created by a single dislocation}\\
Remark that when there is no dislocations, the energy associated to a continuous displacement $U(X)$ for $X=(X_1,X_2)$ is formally
$$E=\int_{\R^2} \frac12 |\nabla U|^2 .$$
Therefore the Euler-Lagrange equation (which is the corresponding equation of elasticity for this model) is
\begin{equation}\label{eq::30}
\Delta U=0 .
\end{equation}
Let us now consider the following function
$$U_0(X)=\frac{1}{2\pi} \arctan \left(\frac{X_1}{X_2}\right) + \frac{1}{4} \mbox{sgn } (X_2)$$
where $\mbox{sgn}$ is the sign function.
This function satisfies
$$\left\{\begin{array}{l}
U_0(X_1,X_2)=-U_0(X_1,-X_2)\\
\\
U_0(X_1,0^+)=-U_0(X_1,0^-) \to 0 \quad \mbox{as}\quad X_1\to -\infty\\
\\
\displaystyle{U_0(X_1,0^+)=-U_0(X_1,0^-) \to \frac12 \quad \mbox{as}\quad X_1\to +\infty .}
\end{array}\right.$$
Moreover we can easily check that
$$\mbox{div }\left(\nabla U_0 - H(X_1)\delta_0(X_2)e_2\right)=0 \quad \mbox{in}\quad {\mathcal D}'(\R^2)$$
where $H$ is the Heavyside function and $\delta_0$ is the Dirac mass.
This equation is the analogue of equation (\ref{eq::30}) when there is a dislocation.
This shows that in a continuum mechanics framework associated to the particular lattice that we consider, 
the function $U_0$ is the displacement corresponding to a dislocation with Burgers vector $e_1$.
In particular, the stress created by this dislocation is then given by
$$\sigma^0=\frac{1}{2\pi}\left(\begin{array}{cc}
\displaystyle{\frac{X_2}{X_1^2+X_2^2}} & \displaystyle{-\frac{X_1}{X_1^2+X_2^2}}\\
\\
\displaystyle{-\frac{X_1}{X_1^2+X_2^2}} & 0
\end{array}\right)$$
and then
\begin{equation}\label{eq::31}
\sigma_{12}^0(X_1,0)=-\frac{1}{2\pi X_1}
\end{equation}
which is the asymptotic shear stress at the point $(X_1,0)\in\R^2$
created by a single dislocation positioned at the origin, and with Burgers vector $e_1$.

%%%%%%%%%%%%%%%%%%%%%%%%%%%%%%%%%%%%%%%%%%%%%%%%%%%%%%%%%%%%%%%%%%%%%%%%%%
\subsection{Rescaling of the generalised FK model}\label{eps1_1}

In this subsection, in order to simplify the notation we denote by $\varepsilon:=\e_{1}>0$ 
the small parameter in the first passage of the scheme
(\ref{tableau}).  
We are interested in the case of asymptotically small potential
$\widetilde{W}$ for which we expect an asymptotically 
large dislocation core. This means that in this limit, we expect to be
able to describe the discrete displacement $U_I$ 
by a continuous function.

More precisely, we first define 
the rescaled integer coordinates:
$$\Omega^\varepsilon= \left(\varepsilon \Z\right) \times \varepsilon \left(\N\backslash \left\{0\right\}\right),\quad \partial \Omega^\varepsilon = \left(\varepsilon \Z\right) \times \left\{0\right\} .$$
%, \quad \overline{\Omega}^\varepsilon= \Omega^\varepsilon \cup \partial \Omega^\varepsilon$$
Then we write the potential as
$$\widetilde{W}=\frac{\varepsilon}{2} W^\varepsilon$$
and define the rescaled function
$$u^\varepsilon(X,t)=2U_{\frac{X}{\varepsilon}}\left(\frac{t}{\varepsilon}\right) \quad \mbox{for}\quad X=(X_1,X_2)\in \overline{\Omega}^\varepsilon,\quad  t\in[0,+\infty) .$$
Remark that the factor $2$ in the definition of $u^\varepsilon$ permits to interprete $u^\varepsilon$ as the jump of the  displacement in the direction $e_1$, when we pass from hyperplane $X_2=-\varepsilon$ to the hyperplane $X_2=0$.

We can easily check that $u^\varepsilon$ solves the following system of
equations (with the particular value $\sigma=0$)

\begin{equation}\label{eq::3}
\left\{\begin{array}{l}
\displaystyle{0=\frac{1}{\varepsilon^2}\sum_{J\in\Z^2, \ |J|=1} (u^\varepsilon(X+\varepsilon J,t)-u^\varepsilon(X,t))} \ \quad \mbox{for all}\quad (X,t)\in \Omega^\varepsilon \times (0,+\infty)\\
\\
u^\varepsilon_t(X,t)=\e\e_{2}{\sigma}\left(\e\e_{2}X_1\right)-(W^\varepsilon)'(u^\varepsilon(X,t)) +I^\varepsilon[u^\varepsilon](X,t)  \quad \mbox{for all}\quad (X,t)\in \partial\Omega^\varepsilon \times (0,+\infty)\\
\\
\displaystyle{\mbox{with}\quad I^\varepsilon[u^\varepsilon](X,t)= \frac{1}{\varepsilon}\sum_{J\in\Z^2,\ |J|=1,\  J_2\ge 0} (u^\varepsilon(X+\varepsilon J,t)-u^\varepsilon(X,t))}.
\end{array}\right.
\end{equation}
Here $\e_{2}>0$ is a small parameter, and the scalar function ${\sigma}$ has been introduced  
to take into account the possible external applied  shear stress on the material.
We will also assume that the initial data satisfies
\begin{equation}\label{eq::4}
u^\varepsilon(X,0)=u_0(X) \quad \mbox{for all}\quad X\in \partial \Omega^\varepsilon   
\end{equation} 
where $u_0$ is a given function independent on $\varepsilon$ and smooth enough.

In order to identify a limit model as $\varepsilon$ goes to zero, we also make the following assumption
\begin{equation}\label{eq::2}
||W^\varepsilon -W||_{C^1(\R)} \to 0 \quad \mbox{as}\quad \varepsilon \to 0
\end{equation}
for some new potential $W$ satisfying the following assumption:\\
\noindent {\bf Assumption (A1)}
$$\left\{\begin{array}{l}
\mbox{The potential} \quad $W$\quad  \mbox{satisfies} \quad \mbox{\bf ($\widetilde{A1}$) i), ii), iii)}\\
\mbox{\bf iv) (Non degenerate minima):}\quad \alpha:=W''(0)>0 .
\end{array}\right.$$

In (\ref{eq::2}), we use the $C^1$ norm, 
because this is the first derivative of the potential that appears in the equations.
Remark that condition (\ref{eq::2}) can be fulfilled, if we assume for
instance that $\widetilde{W}$ satisfies assumption
$(\widetilde{\mbox{A1}})$ with $\delta=\delta_\varepsilon < <
\varepsilon$. We also make the following assumption on the stress:\\
\noindent {\bf Assumption (A2)}\\
There exists a constant $C>0$ such that
$$|\sigma|+|\sigma_x|+|\sigma_{xx}|\le C \quad \mbox{for all}\quad x\in\R .$$

%%%%%%%%%%%%%%%%%%%%%%%%%%%%%%%%%%%%%%%%%%%%%%%%%%%%%%%%%%%%%%%%%%%%%%%%%%
%%%%%%%%%%%%%%%%%%%%%%%%%%%%%%%%%%%%%%%%%%%%%%%%%%%%%%%%%%%%%%%%%%%%%%%%%%
\section{The Peierls-Nabarro model}

%%%%%%%%%%%%%%%%%%%%%%%%%%%%%%%%%%%%%%%%%%%%%%%%%%%%%%%%%%%%%%%%%%%%%%%%%%
\subsection{Description of the PN model}\label{eps1_2}

In this section we introduce the Peierls-Nabarro model, 
which is a phase field model (see \cite{HL} for a presentation of this model). In this model, 
phase transitions describe the dislocation cores. 
We set
$$\Omega=\left\{X=(X_1,X_2)\in\R^2, \quad X_2 >0\right\} .$$
A function $u^0(X,t)$ is said to be a solution of the PN model, if it satisfies the following system
\begin{equation}\label{eq::5}
\left\{\begin{array}{ll}
\displaystyle{0= \Delta u^0} &\quad \mbox{on}\quad  \Omega \times (0,+\infty)\\
\\
\displaystyle{u^0_t=2\e_{2}{\sigma}\left(\e_{2}X_1\right)-W'(u^0) + \frac{\partial u^0}{\partial X_2}}  &\quad \mbox{on}\quad \partial\Omega \times (0,+\infty)
\end{array}\right.
\end{equation}
with initial data
\begin{equation}\label{eq::6}
u^0(X,0)=u_0(X) \quad \mbox{for all}\quad X\in \partial\Omega .
\end{equation}

\begin{figure}[h]
\centering\epsfig{figure=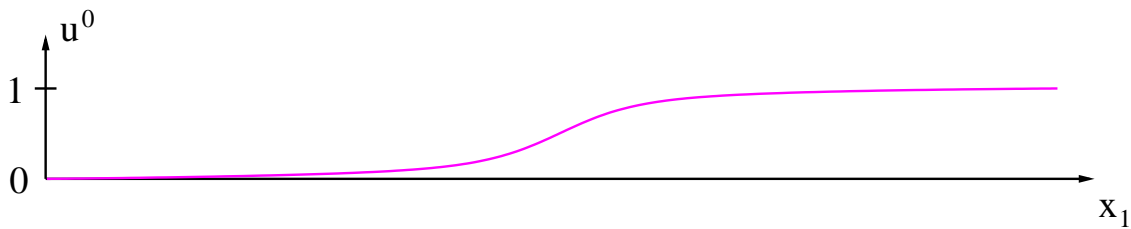,width=150mm}
\caption{Phase transition for a edge dislocation with Burgers vector $e_1$ for $X_2=0$}\label{F3}
\end{figure}

The stationary version of this model has been originally introduced 
in order to propose a method to compute at the equilibrium a finite stress created by a dislocation.
In this model, $u^0$ is the phase transition. For instance, for a edge dislocation with Burgers vector $e_1$ as presented in Section \ref{s1}, $u^0$ is a transition between the value $0$ on the left to the value $1$ on the right (see Figure \ref{F3}). In the special case $u^0_t=0={\sigma}$ and for sinusoidal potentials $W$, 
the stationary solution $u^0$ is known explicitely (see for instance
\cite{CSM}), which makes the PN model very attractive. Let us mention
that a physical and numerical study of the evolution problem
(\ref{eq::5}) has been treated in \cite{MBW1998}.

\begin{rem}
Remark that when we consider system (\ref{eq::5})-(\ref{eq::6}) in the framework of viscosity solutions,
the evolution equation on the boundary $\partial\Omega$ appears to be a boundary condition of the system.
For this reason, as it is usual for viscosity solutions (see for instance \cite{B0,BDL}), this boundary condition has to be understood technically 
in the sense that on $\partial \Omega$ the function $u^0$ solves pointwisely either $0=\Delta u^0$ or the evolution equation.
\end{rem}

%%%%%%%%%%%%%%%%%%%%%%%%%%%%%%%%%%%%%%%%%%%%%%%%%%%%%%%%%%%%%%%%%%%%%%%%%%
\subsection{Convergence of the generalised FK  model to the PN model}\label{eps1_3}

We have the following result
\begin{theorem}{\bf (Formal convergence of FK to PN, $\e=\e_1 \to 0$)}\label{th::1}\\
Let $\varepsilon>0$.  For the initial data $u_0\in W^{2,\infty}(\overline{\Omega})$ which is assumed harmonic on $\Omega$,
and under assumption ($\widetilde{A1}$) on $\varepsilon W^\varepsilon$,
and (A2) on $\sigma$, there exists a unique viscosity solution $u^\varepsilon$ 
of system (\ref{eq::3})-(\ref{eq::4}). Moreover assuming (\ref{eq::2})
with the potential $W$ satisfying assumption (A1),
 then, as $\varepsilon$ goes to zero, the solution $u^\varepsilon$
 formally converges  to a viscosity solution of system 
(\ref{eq::5})-(\ref{eq::6}). 
\end{theorem}

The proof of Theorem \ref{th::1} is done in full details in \cite{FinoIM}.\\

\noindent {\bf Sketch of the proof of Theorem \ref{th::1}}\\
One way to guess the limit model (\ref{eq::5})-(\ref{eq::6}) is to pass to the limit formally in system 
(\ref{eq::3})-(\ref{eq::4}) assuming that the solution $u^\varepsilon$ (and its derivatives) converges to a limit $u^0$.
The convergence in the system is then obtained using a simple Taylor expansion.
The existence of a solution $u^\varepsilon$ to system (\ref{eq::3})-(\ref{eq::4}) is technically delicate 
and is based on the proof of a suitable comparison principle for this system.

%%%%%%%%%%%%%%%%%%%%%%%%%%%%%%%%%%%%%%%%%%%%%%%%%%%%%%%%%%%%%%%%%%%%%%%%%%
\subsection{Reformulation of the PN model}
In this subsection and in Subsection \ref{eps2_2}, 
in order to simplify the notation we denote by $\varepsilon:=\e_{2}>0$ 
the small parameter in the second passage of the scheme
(\ref{tableau}).  We recall that it is well known that for bounded
smooth  functions $u^0$ 
defined on $\overline{\Omega}$ which are harmonic on $\Omega$, we can write
$$\frac{\partial u^0}{\partial X_2}(X_1,0) = L(u^0(\cdot,0))(X_1)\quad
\mbox{for all}\quad (X_1,0)\in\partial
 \Omega$$
where for a general  bounded smooth function $w$, the linear operator
$L$ is given by the Levy-Khintchine formula 
(see Theorem 1 in \cite{DI}):
\begin{equation}\label{eq::17}
(Lw)(x)=\frac{1}{\pi}\int_\R \frac{dz}{z^2}\left(w(x+z)-w(x)-zw'(x) 1_{\left\{|z|\le 1\right\}}\right) .
\end{equation}

Then for smooth solutions $u^0$, system (\ref{eq::5}) can be rewritten
for $V(X_1,t)=u^0(X_1,X_2,t)_{|X_2=0}$  with $x=X_1\in\R$ as
\begin{equation}\label{eq::7}
V_t = 2\e\sigma(\e x)-W'(V) + LV \quad \mbox{on}\quad \R.
\end{equation}
We also recall (see \cite{CSM}) that there exists a unique function $\phi$ solution of 
\begin{equation}\label{eq::12}
\left\{\begin{array}{l}
0=L\phi -W'(\phi) \quad \mbox{on}\quad \R\\
\\
\displaystyle{\phi'>0 \quad \mbox{and}\quad \phi(-\infty )=0,\quad \phi(0)=\frac12,\quad \phi(+\infty)=1} .
\end{array}\right.
\end{equation}

The function $\phi$ is called the layer solution and a translation of
$\phi$ is pictured on Figure \ref{F3}.

%%%%%%%%%%%%%%%%%%%%%%%%%%%%%%%%%%%%%%%%%%%%%%%%%%%%%%%%%%%%%%%%%%%%%%%%%%
\subsection{Rescaling of the PN model}\label{eps2_2}

We now consider the following rescaling
$$v^\varepsilon(x,t)=V\left(\frac{x}{\varepsilon},\frac{t}{\varepsilon^2}\right) .$$
Then system (\ref{eq::7}) can be rewritten as
\begin{equation}\label{eq::8}
v^\varepsilon_t = \frac{1}{\varepsilon}\left\{ Lv^\varepsilon -\frac{1}{\varepsilon}W'(v^\varepsilon) + 2\sigma(x)\right\} \quad \mbox{on}\quad \R
\end{equation}
with initial condition
\begin{equation}\label{eq::9}
v^\varepsilon(x,0)=v_0^\varepsilon(x) \quad \mbox{for}\quad x\in\R .
\end{equation}

Again, a good notion of solution for system (\ref{eq::8})-(\ref{eq::9}) 
is the notion of viscosity solution for non local equations (see for instance \cite{BI}).

Here we will choose carefully the initial condition $v^\varepsilon_0$ as follows\\
\noindent {\bf Assumption (A3)}
$$\left\{\begin{array}{l}
x^0_1< x_2^0< ... < x_N^0\\
\\
\displaystyle{v^\varepsilon_0(x)=\frac{\varepsilon}{\alpha}\cdot 2\sigma(x)+\sum_{i=1}^N \phi\left(\frac{x-x_i^0}{\varepsilon}\right)}
\end{array}\right.$$
where we recall that $\alpha=W''(0)>0$ and $\phi$ is defined in (\ref{eq::12}).

%%%%%%%%%%%%%%%%%%%%%%%%%%%%%%%%%%%%%%%%%%%%%%%%%%%%%%%%%%%%%%%%%%%%%%%%%%
%%%%%%%%%%%%%%%%%%%%%%%%%%%%%%%%%%%%%%%%%%%%%%%%%%%%%%%%%%%%%%%%%%%%%%%%%%
\section{Dynamics of discrete dislocations}

%%%%%%%%%%%%%%%%%%%%%%%%%%%%%%%%%%%%%%%%%%%%%%%%%%%%%%%%%%%%%%%%%%%%%%%%%%
\subsection{Description of the DDD model}\label{eps2_3}

\begin{figure}[h]
\centering\epsfig{figure=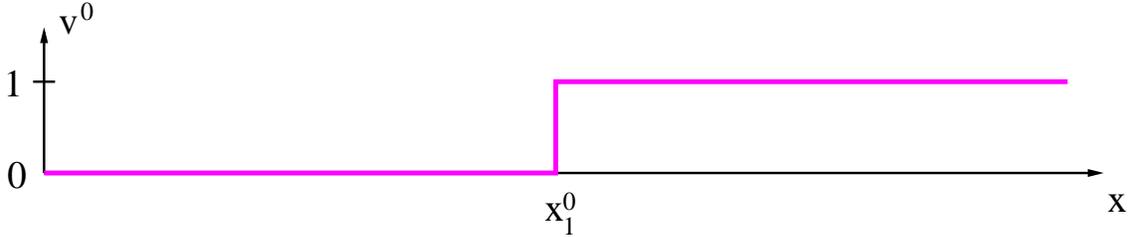,width=150mm}
\caption{Sharp interface describing a discrete dislocation at $x=x_1^0$}\label{F4}
\end{figure}

In this section we assume that the phase transition reduces to a sharp interface where the transition is localized at the position $x=x^0_1 \in\R$. For a dislocation associated to a Burgers vector $e_1$, the sharp interface is associated to a non-decreasing step function like the one of Figure \ref{F4}.

More generally, we can consider the case of $N$ dislocations (or particules)
of positions $(x_i(t))_{i=1,...,N}$ solving the following system
\begin{equation}\label{eq::10}
\frac{dx_i}{dt}=-\gamma\left(\sigma(x_i) 
+\sum_{j\not= i} V'(x_i-x_j)\right) \quad \mbox{on}\quad (0,+\infty) \quad \mbox{for}\quad i=1,...,N
\end{equation}
with the two-body interaction potential
$$\displaystyle{V(x)=-\frac{1}{2\pi} \ln |x|}$$
with initial data
\begin{equation}\label{eq::11}
x_i(0)=x_i^0 \quad \mbox{for}\quad i=1,...,N .
\end{equation}
Here the constant $\gamma>0$ is the inverse of the damping factor. It is 
related to the layer solution $\phi$ defined in (\ref{eq::12}) and is given by
$$\gamma=2\left(\int_{\R}(\phi')^2\right)^{-1} .$$

The function $\sigma$ is the applied shear stress and $V'(x-x_j)$ is the (singular) shear stress created at the point $x$ by the dislocation $x_j$. This corresponds exactly to the shear stress already computed in (\ref{eq::31}).
The total stress $\displaystyle{\sigma(x_i)+\sum_{j\not= i} V'(x_i-x_j)}$ is called the resolved Peach-Koehler force acting on the dislocation $x_i$.

%%%%%%%%%%%%%%%%%%%%%%%%%%%%%%%%%%%%%%%%%%%%%%%%%%%%%%%%%%%%%%%%%%%%%%%%%%
\subsection{Convergence of the PN  model to the DDD model}\label{eps2_4}

We have
\begin{theorem}{\bf (Convergence of PN  to DDD, $\e=\e_2 \to 0$)}\label{th::2}\\
Let $\varepsilon>0$. Under assumptions (A1)-(A2)-(A3), there exists a unique viscosity solution $v^\varepsilon$ of (\ref{eq::8})-(\ref{eq::9}). Moreover there exists a unique solution of (\ref{eq::10})-(\ref{eq::11}), and we define
$$v^0(x,t)=\sum_{i=1,...,N} H(x-x_i(t))$$
where $H$ is the Heavyside function.
Then as $\varepsilon$ goes to zero, the function $v^\varepsilon$ converges to $v^0$ in the following sense
$$\limsup_{(x',t')\to (x,t),\ \varepsilon\to 0} v^\varepsilon(x',t') \le (v^0)^*(x,t)$$
and
$$\liminf_{(x',t')\to (x,t),\ \varepsilon\to 0} v^\varepsilon(x',t') \ge (v^0)_*(x,t) .$$
\end{theorem}

The proof of this result is done in full details in \cite{GM}.

\begin{rem}
We recall that the semi-continuous envelopes of a function $v$ are defined as
$$v^*(x,t)= \limsup_{(x',t')\to (x,t)} v(x',t') \quad \mbox{and}\quad v_*(x,t)= \liminf_{(x',t')\to (x,t)} v(x',t') .$$
\end{rem}

\noindent {\bf Sketch of the proof of convergence}\\
The existence of a solution for all time of the ODE system (\ref{eq::10})-(\ref{eq::11}) comes from the fact 
that $V(x)$ is a convex potential outside the origin. 
This property allows to show that the minimal distance between particles
$$d(t)=\inf_{i\not= j} |x_i(t)-x_j(t)|$$
satisfies
\begin{equation}\label{eq::20}
d(t)\ge d(0)e^{-C\gamma t}
\end{equation}
which prevents the meeting of the particles at any finite time.\\
Then the main idea to prove the convergence is to approximate the solution $v^\varepsilon$ by the following ansatz
$$\displaystyle{\tilde{v}^\varepsilon(x,t)
=\frac{\varepsilon}{\alpha}\cdot 2\sigma(x)+\sum_{i=1}^N \left\{\phi\left(\frac{x-x_i}{\varepsilon}\right)-\varepsilon\dot{x}_i(t)\psi\left(\frac{x-x_i}{\varepsilon}\right) \right\}\quad \mbox{with}\quad \dot{x}_i(t)= \frac{dx_i}{dt}(t)} $$
where $\alpha=W''(0)$ and the corrector $\psi$ solves the following equation
$$L\psi-W''(\phi)\psi =\phi' +\eta \left(W''(\phi)-W''(0)\right) 
\quad \mbox{with}\quad \eta=\displaystyle{\frac{1}{W''(0)}\int_{\R}(\phi')^2} .$$
The stress created in $x$ by a dislocation positioned at the origin, comes from the following property
$$\phi(x)-H(x) \sim -\frac{1}{\alpha\pi x} \quad \mbox{as}\quad |x|\to +\infty .$$
The rest of the proof of  convergence of $v^\varepsilon$ is done by construction of sub and super solutions based 
on the ansatz $\tilde{v}^\varepsilon$.

%%%%%%%%%%%%%%%%%%%%%%%%%%%%%%%%%%%%%%%%%%%%%%%%%%%%%%%%%%%%%%%%%%%%%%%%%%
\subsection{Rescaling of the DDD model}
In this subsection, 
in order to simplify the notation we denote by $\varepsilon:=\e_{3}>0$ 
the small parameter in the  third passage of the scheme
(\ref{tableau}). We consider a given initial data $w_0$ which satisfies\\
\noindent {\bf Assumption (A4)}\\
$$\left\{\begin{array}{l}
w_0\in W^{2,\infty}(\R),\\
w_0'> 0, \quad w_0(-\infty)=0 .
\end{array}\right.$$
We also introduce the integer $N_\varepsilon$ and the position of the dislocations $x_1^0<...<x_{N_\varepsilon}^0$ such that
$$\displaystyle{\sum_{i=1,...,N_\varepsilon} H(x-x_i^0) = \big\lfloor \frac{w_0(\varepsilon x)}{\varepsilon}\big\rfloor}$$
where $\lfloor \cdot \rfloor$ denotes the floor function. We also assume
that the stress $\sigma$ is periodic. Precisely, we make the following assumption:\\
\noindent {\bf Assumption (A2')}
$$\sigma \in C^2(\R)\quad \mbox{and}\quad \sigma(x+1)=\sigma(x) \quad \mbox{for all}\quad x\in\R .$$
This assumption allows
to study the collective behaviour of dislocations in a landscape with periodic obstacles, 
and to get the effective macroscopic model by a periodic homogenization approach.

Then we consider the solution $(x_i(t))_{i=1,...,N_\varepsilon}$ of 
the system (\ref{eq::10})-(\ref{eq::11}) with $N=N_\varepsilon$ and define the function
$$v^0(x,t)=\sum_{i=1,...,N_\varepsilon} H(x-x_i(t))$$
and the rescaling
\begin{equation}\label{eq::13}
w^\varepsilon(x,t)=\varepsilon v^0\left(\frac{x}{\varepsilon},\frac{t}{\varepsilon}\right) .
\end{equation}

%%%%%%%%%%%%%%%%%%%%%%%%%%%%%%%%%%%%%%%%%%%%%%%%%%%%%%%%%%%%%%%%%%%%%%%%%%
%%%%%%%%%%%%%%%%%%%%%%%%%%%%%%%%%%%%%%%%%%%%%%%%%%%%%%%%%%%%%%%%%%%%%%%%%%
\section{Dislocation density model}

%%%%%%%%%%%%%%%%%%%%%%%%%%%%%%%%%%%%%%%%%%%%%%%%%%%%%%%%%%%%%%%%%%%%%%%%%%
\subsection{Description of the DD model}

We first introduce a function $g:(0,+\infty)\times \R\to \R$ which satisfies\\
\noindent {\bf Assumption (A5)}
$$\left\{\begin{array}{l}
g\in C^0((0,+\infty)\times \R),\\
\\
l\mapsto g(\rho,l) \quad \mbox{is nondecreasing.}
\end{array}\right.$$
Then we consider a function $w^0(x,t)$ which is a solution of 
\begin{equation}\label{eq::15}
w^0_t=g(w^0_x,Lw^0) \quad \mbox{on}\quad \R\times (0,+\infty)
\end{equation} 
where the operator $L$ is defined in (\ref{eq::17}), and 
with initial data
\begin{equation}\label{eq::16}
w^0(x,0)=w_0(x) \quad \mbox{for all}\quad x\in\R .
\end{equation}

Here the function $w^0$ is such that its derivative $w^0_x$ represents the macroscopic dislocation density.
Moreover $w^0$ can be seen as the plastic strain localized in plane
$x_2=0$ and $\displaystyle{\frac12 Lw^0}$
 can be identified to the stress created by the dislocation density $w^0_x$. 
Equation (\ref{eq::15}) can be interpreted as the plastic flow rule in a
model for macroscopic crystal plasticity. Indeed, from a mechanical
point of view, we have the following table (see also \cite{IMR}) of equivalence between our
homogenized model and a classical model in mechanics for
elasto-visco-plasticity of crystals (see \cite{FPZ}).

\begin{center}
\noindent \begin{tabular}{||l|r|r||}\hline
&  %Mechanics
Crystal elasto-visco-plasticity &  DD model \\ \hline\hline
\begin{tabular}{l} resolved  plastic strain \end{tabular}&
$\gamma(x_1)\delta_0(x_2)$ &  $w^0(x_1)\delta_0(x_2)$\\\hline
\begin{tabular}{l} Nye  tensor of\\ dislocations densities
\end{tabular}
& $\alpha=(e_1\otimes e_2)  \gamma'(x_1) \delta_0(x_2)$ &
$\alpha=(e_1\otimes e_2) w^0_x(x_1) \delta_0(x_2)$ \\\hline
\begin{tabular}{l} exterior applied stress\end{tabular}  & $\Sigma^{ext}$ &   \\\hline
\begin{tabular}{l} microscopic resolved \\ shear stress 
\end{tabular}
& & $\sigma-\displaystyle{\int_{(0,1)}
\sigma}$ \\ \hline
\begin{tabular}{l} resolved  exterior \\ applied stress  
\end{tabular} & $\Sigma^{ext}: e^0$ & $\displaystyle{\int_{(0,1)}
\sigma}$ \\ \hline
%   & &  \\
\begin{tabular}{l} displacement  \end{tabular}  & $v=v_1e_1$ &  \\\hline
\begin{tabular}{l} strain \end{tabular}  & $\displaystyle{e:=e(v)-e^0 \gamma \delta_0(x_2)}$ &  \\
  & with $\displaystyle{e(v) := \frac12\left(\nabla v+{}^t\nabla v\right)}$ &  \\
 & and $\displaystyle{e^0:=\frac12 (e_2\otimes e_1 + e_1\otimes e_2)}$ &  \\\hline
\begin{tabular}{l} total elastic energy \end{tabular}  &
$\displaystyle{E:=\int_{\R^2} \frac12 (\Lambda :e):e + \Sigma^{ext}: e}$
&  $\displaystyle{E:=\int_{\R} -\frac14 w^0  L w^0 - \left(\int_{(0,1)}\sigma\right)w^0}$ \\\hline
\begin{tabular}{l} macroscopic stress  \end{tabular}  & $\Sigma:= \Lambda : e + \Sigma^{ext}$ &  \\\hline
\begin{tabular}{l} resolved macroscopic\\ shear stress  
\end{tabular}
& $\tau := \Sigma : e^0$ & $\tau:=\displaystyle{ \frac12 L w^0 + \int_{(0,1)} \sigma}$   \\\hline
\begin{tabular}{l} visco-plastic law  \end{tabular} & $\displaystyle{\frac{\partial \gamma}{\partial t}=f(\tau)}$ &
$\displaystyle{\frac{\partial w^0}{\partial t}=g\left( w^0_x,
  Lw^0 + \int_{(0,1)} 2\sigma\right)}$\\\hline
\begin{tabular}{l} energy decay \end{tabular}& $\displaystyle{\frac{d}{dt}E=\int_{\R} -\tau f(\tau) \le
 0}$ & $\displaystyle{\frac{d}{dt}E=\int_{\R} - \tau
   g\left( w^0_x,2\tau \right)}\le 0$  \\\hline
%    & &  \\\hline
%    & &  \\\hline
%\hline
\end{tabular}\hfill\break\hfill\break
\end{center}

\begin{figure}[!h]
\centering\epsfig{figure=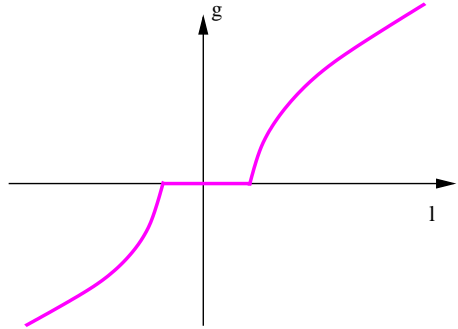,width=60mm}
\caption{Graph of the map $l\mapsto g(\rho,l)$ showing a threshold effect}\label{F5}
\end{figure}

Remark that when we choose the microscopic stress $\sigma$ in Assumption
(A2') so that $\displaystyle{\int_{(0,1)}\sigma=0}$, we expect that $g$
satisfies 
\begin{equation}\label{g=0}
g(\rho,0)=0.
\end{equation}
This equality reflects the pinning of dislocations (see \cite[Th
2.6]{FIM}). In the model presented in the previous table, the plastic strain
velocity $w^0_t$ is prescribed by the function $g$ 
(whose typical profile can be seen on Figure \ref{F5}) which is assumed
to satisfy (\ref{g=0}).

Let us mention that a similar profile  for $g$ as a
function of $l$ only have been obtained rigorously for a different model
involving the motion of a phase transition solution of Allen-Cahn type
equations (see \cite{DY2006}). 

%%%%%%%%%%%%%%%%%%%%%%%%%%%%%%%%%%%%%%%%%%%%%%%%%%%%%%%%%%%%%%%%%%%%%%%%%%
\subsection{Convergence of the DDD  model to the DD model}

We have the following result
\begin{theorem}{\bf (Convergence of DDD  to DD, $\e=\e_3 \to 0$)}\label{th::3}\\
Let us assume (A2'), (A4). Then there exists a function $g$ satisfying assumption (A5).
Moreover the function $w^\varepsilon$ defined in (\ref{eq::13}) converges to the unique solution $w^0$ of (\ref{eq::15})-(\ref{eq::16}), locally uniformly on $\R\times [0,+\infty)$.
\end{theorem}

The proof of this result is done in full details in \cite{FIM}.

Remark that Theorem \ref{th::3} is an homogenization result in the periodic setting. In the particular case where the periodic stress $\sigma$ is equal to zero, we get
$$g(\rho,l)=\frac{\gamma}{2} \rho l .$$
The presence of a non-zero $1$-periodic stress with zero mean value, creates a threshold phenomenon 
where for a fixed dislocation density $\rho$, the quantity $g(\rho,l)$ can be equal to zero if $|l|$ is small enough
(see for instance the numerical simulations in \cite{EIM} which look like Figure \ref{F5}).\\

\noindent {\bf Sketch of the proof of Theorem \ref{th::3}}\\
\noindent {\bf Step 1 : Formal determination of the function $g$}\\
To determine the function $g$, we can look formally for  $x_i(t)$ defined for all $i\in\Z$, which are particular solutions of the ODE system (\ref{eq::10}) with $\sigma$ replaced by $\displaystyle{\frac{l}{2}+\sigma}$, such that
$$x_i(t)=h(vt + i/\rho) \quad \mbox{with}\quad h(a+1)=1+h(a) \quad \mbox{for all}\quad a\in\R$$
where such a function $h$ is called a hull function. Both $h$ and the constant $v$ have to be determined.
It can be shown that $v$ is unique. Then we set
$$\displaystyle{g(\rho,l)= - v \rho}$$
which is known in physics as the Orowan's law.\\
\noindent {\bf Step 2: Regularization at short distances}\\
To avoid the singularity of the potential $\displaystyle{V(x)=-\frac{1}{2\pi}\ln |x|}$, 
we can first approximate it by the following  symmetric and continuous potential
$$V_\delta(x)=\left\{\begin{array}{l}
V(x) \quad \mbox{if}\quad |x|\ge \delta\\
\mbox{linear if}\quad x\in (-\delta,\delta)\backslash \left\{0\right\} .
\end{array}\right.$$
We consider the function $v^0(x,t)=\sum_{i=1,...,N_\varepsilon} H(x-x_i(t))$ 
associated to the dynamics (\ref{eq::10}) 
where the potential $V$ is replaced by $V_\delta$.
Then it is possible to show that $v^0$ satisfies the following equation
with $l=0$
\begin{equation}\label{eq::21}
\displaystyle{v^0_t = |v^0_x| \gamma \left\{\sigma(x)+\frac{l}{2} + M_\delta[v^0(\cdot,t)](x)\right\}}
\end{equation}
where for a general function $w(x)$, we can define the non local operator
$$M_\delta[{w}](x)=\int_{|z|> \delta} dz\ V_\delta''(z) E\left({w}(x+z)-{w}(x)\right) \quad \mbox{with}\quad E(a)=\frac12 + k \quad \mbox{if}\quad k\le a < k+1,\quad k\in\Z$$
where $E$ is a odd integer part function.
This is possible to introduce a suitable good notion of viscosity solution for equation (\ref{eq::21}) (see \cite{FIM}).
In particular, we can show that if $v^0(x,0)=\rho x$, then $v^0(x,t)/t \to g_\delta(\rho,l)$ as $t\to +\infty$.
Moreover it is possible to show the following estimate
\begin{equation}\label{eq::22}
|g_\delta(\rho,l)-g(\rho,l)|\le \frac{C(\rho)}{|\ln \delta|} .
\end{equation}
\noindent {\bf Step 3: Sketch of the proof of convergence in the regularized case}\\
After a rescaling of the solution $v^0$ of (\ref{eq::21}) with $l=0$, we see that $w^\varepsilon(x,t)=\varepsilon v^0(x/\varepsilon,t/\varepsilon)$ solves an equation
$$\displaystyle{w^\varepsilon_t = |w^\varepsilon_x| \gamma \left\{\sigma(x/\varepsilon) + M_\delta^\varepsilon[w^\varepsilon(\cdot,t)](x)\right\}}$$
for some rescaled non local operator $M_\delta^\varepsilon$. More generally, any continuous solution $w^\varepsilon$ of the previous equation, can be formally written as
$$w^\varepsilon(x,t)\simeq w^0(x,t)+\varepsilon r(x/\varepsilon)$$
where $r$ is a suitable corrector. One fundamental remark is that as $\varepsilon$ goes to zero, 
we can asymptotically split the non local term 
$$M_\delta^\varepsilon[w^0(\cdot,t) + \varepsilon r(\cdot /\varepsilon)](x)\quad \simeq \quad \frac12 (Lw^0)(x) \quad + \quad S[r,w^0_x(x)](x/\varepsilon)$$
into its long range contribution  $\frac12 Lw^0$ and a short range contribution $S$ involving the corrector $r$.
Remark that this long range contribution $\frac12 Lw^0 = l/2$ is related to the introduction 
of the constant $l/2$ into equation (\ref{eq::21}) used in the definition of $g_\delta(\rho,l)$.
Taking into account this asymptotical splitting, it is then possible to show the convergence of $w^\varepsilon$
to the solution of (\ref{eq::15})-(\ref{eq::16}) with $g$ replaced by $g_\delta$. The proof can be done in the framework of viscosity solutions, adapting the Evans' perturbed test function method.\\
\noindent {\bf Step 4: Sketch of the proof of convergence in the singular case}\\
The singular case can be reached using an approximation argument. 
On the one hand, estimate (\ref{eq::20}) insures that the dynamics (\ref{eq::10}) on the time interval $(0,T/\varepsilon)$
is equivalent to the same dynamics with $V$ replaced by $V_\delta$ for $\delta \le \delta_\varepsilon= d(0)e^{-\gamma CT/\varepsilon}$. On the other hand,
estimate (\ref{eq::22}) is independent on $\varepsilon$. Then choosing $\delta=\delta_\varepsilon$, the  convergence of the solution $w^\varepsilon$ on the time interval $(0,T)$ can then be obtained by an adaptation of the arguments in the regularized case.

%%%%%%%%%%%%%%%%%%%%%%%%%%%%%%%%%%%%%%%%%%%%%%%%%%%%%%%%%%%%%%%%%%%%%%%%%%
\section{Conclusion}

We considered a two-dimensional Frenkel-Kontorova model in the fully overdamped regime.
From this model, we derived by a scaling argument the time-dependent Peierls-Nabarro model.
Looking at the sharp interface limit of the phase transitions of the Peierls-Nabarro model,
we were able to identify a dynamics of particles that corresponds to the classical discrete dislocation dynamics,
in the particular case of parallel straight edge dislocation lines in the same glide plane with the same Burgers vector.
Considering the motion of these particles in a landscape with periodic obstacles, we were able to identify at large scale
an evolution model for the dynamics of a density of dislocations. This model is a macroscopic model for 
crystal visco-elasto-plasticity, where we predicted a plastic flow rule.
 This last model shows in particular a threshold effect where
 dislocations can be pinned in the obstacles, if the effective stress
 acting on these dislocations is too small.

In order to present a summary of our approach, we give here a  diagram
(see Figure 6.1) that
shows the links between the four models treated in this paper.

Up to our knowledge, this derivation of classical models 
from a single microscopic model (the 2D Frenkel-Kontorova model),
seems new and allows to make clear connections between different
modelling of dislocation dynamics.

\newpage
$$\fbox{ 
$\begin{array}{ll} \hspace{5cm} \mbox{2D Frenkel-Kontorova (FK)}\\\\
\left \{\begin{array}{ll}
\displaystyle{0=\sum_{J \in \Z^2,\; |J|=1} (U_{X+J}(t) -U_{X}(t))} &
\mbox{ for}\quad X=
(X_1,X_2) \quad\mbox{with}\quad X_2\in  \N\backslash \left\{0\right\}\\\\
\displaystyle{ \frac{d}{dt} U_X(t)=  \e_1 \e_2 {\sigma}(\e_1
  \e_2 X_1) - \frac {\e_1}{2} (W^{\e_1})' (2U_X(t))+ 
  I^{1}[U_X(t)]} & \mbox{ for}\quad X_2 =0\\\\
 \mbox{where}\quad I^{1}[U_X(t)]=\displaystyle{\sum_{J \in \Z^2, 
|J|=1,\; J_2 \ge 0} U_{X+J}(t) -U_{X}(t)}
\end{array}\right.
\end{array}
$ }$$

$$\hspace{7cm}\begin{array}{llll}
\e_1    \\
\downarrow\\
0
\end{array}
  \left \Downarrow\begin{array}{llll}    
u^{\e_1}(X,t)= 2 U_{\frac {X}{\e_1}}\left(\frac {t}{\e_1}\right) \to  u^0(X,t), \quad \mbox{as}\quad \e_1 \to 0 \\
  W^{\e_1} \to W, \quad \mbox{as}\quad \e_1 \to 0 \\
  v(X_1,t) = u^0(X_1,X_2,t)_{|X_2=0}; \quad x=X_{1} \\
  \; 
  \end{array} \right . $$

$$ \fbox{ $\begin{array}{ll}\hspace{2cm} \mbox{Peierls-Nabarro (PN)}\\\\
v_t (x,t)=2 \e_2 \sigma (\e_2 x ) - W'(v(x,t)) + Lv(x,t)  
\end{array}
$}$$

$$\hspace{6.5cm}\begin{array}{llll}
\e_2    \\
\downarrow\\
0
\end{array}
  \left \Downarrow\begin{array}{llll} 
  \; \; \\   
v^{\e_2}(x,t)=  v\left(\frac {x}{\e_2} , \frac {t}{\e_2^2}\right) \to  v^0(x,t), \quad \mbox{as}\quad \e_2 \to 0 \\
\displaystyle{ v^0(x,t)= \sum_{i=1,\cdots,N}H(x-x_i(t))}\\
 \;
\end{array} \right . $$

 $$ \fbox{ $\begin{array}{ll} \hspace{1.5cm} \mbox{Dynamics of Discrete
       Dislocations (DDD)}\\\\
 \displaystyle{\dot{x}_i (t)= - \gamma  \left( \sigma (x_i)-\sum_{i
       \neq j} \frac {1}{x_i(t)-x_j(t)}\right)} 
 \quad \mbox{for}\quad i=1, \cdots, N
 \end{array}
 $ }$$
 
$$\hspace{7cm} \begin{array}{llll}
\e_3    \\
\downarrow\\
0
\end{array}
  \left \Downarrow\begin{array}{llll} 
  \; \; \\   
w^{\e_3}(x,t)= \e_3v^0\left(\frac {x}{\e_3} , \frac {t}{\e_3}\right) \to  w^0(x,t), \quad \mbox{as}\quad \e_3 \to 0 \\
 \sigma(x+1)=\sigma(x)\\
 \;
\end{array} \right . $$

$$ \fbox{$\begin{array}{ll} \mbox{Dislocations Density (DD)}\\\\
 w^0_t= g(w^0_x, Lw^0)
 \end{array}
$} $$\\

\noindent{\small {\bf Figure 6:} Descriptive diagram summarizing the
  links between the four models}
\\

%%%%%%%%%%%%%%%%%%%%%%%%%%%%%%%%%%%%%%%%%%%%%%%%%%%%%%%%%%%%%%%%%%%%%%%%%
%%%%%%%%%%%%%%%%%%%%%%%%%%%%%%%%%%%%%%%%%%%%%%%%%%%%%%%%%%%%%%%%%%%%%%%%%
\bigskip
\bigskip

%%%%%%%%%%%%%%%%%%%%%%%%%%%%%%%%%%%%%%%%%%%%%%%%%%%%%%%%%%%%%%%%%%%%%%%%%
%%%%%%%%%%%%%%%%%%%%%%%%%%%%%%%%%%%%%%%%%%%%%%%%%%%%%%%%%%%%%%%%%%%%%%%%%

%%%%%%%%%%%%%%%%%%%%%%%%%%%%%%%%%%%%%%%%%%%%%%%%%%%%%%%%%%%%%%%%%%%%%%%%%
%%%%%%%%%%%%%%%%%%%%%%%%%%%%%%%%%%%%%%%%%%%%%%%%%%%%%%%%%%%%%%%%%%%%%%%%%

\noindent {\bf Aknowledgements}\\
The authors would like to thank L. Truskinovsky for his fruitful remarks
and suggestions which improved the final presentation of the paper and
for interesting references related to this topic. 
This work was supported by the ANR ``Mouvements d'Interfaces, Calcul et Applications''  (2006-2009).

%%%%%%%%%%%%%%%%%%%%%%%%%%%%%%%%%%%%%%%%%%%%%%%%%%%%%%%%%%%%%%%%%%%%%%%%%%%%%%%
%%%%%%%%%%%%%%%%%%%%%%%%%%%%%%%%%%%%%%%%%%%%%%%%%%%%%%%%%%%%%%%%%%%%%%%%%%%%%%%

\end{document}